\documentclass[final]{siamltex}

\newtheorem{remark}{Remark}
\usepackage{epsfig}
\usepackage{amssymb}
\usepackage{amsmath}
\usepackage{showkeys}

\def\Erfc{{\rm Erfc}}

\def\tfrac#1#2{{{\lower.6ex
\hbox{$\scriptstyle#1$}}\over
{\raise.7ex
\hbox{$\scriptstyle#2$}}}}

\def\sign{{\rm sign}}
\def\erfc{{\rm erfc}}
\def\inverfc{{\rm inverfc}}

\def\erf{{\rm erf}}
\def\inverf{{\rm inverf}}

\def\bigO{{\cal O}}

\def\dsp#1{\displaystyle#1}

\def\protectbold#1{\protect{\boldmath{$#1$}}}

\def\Frac#1#2{\frac{\displaystyle{#1}}{\displaystyle{#2}}}

\def\bigO{{\cal O}}

\def\sign{{\rm sign}}

\def\sign{{\rm sign}}
\def\erfc{{\rm erfc}}
\def\erf{{\rm erf}}
\def\bigO{{\cal O}}

\title{Efficient and accurate algorithms for the computation and inversion 
of the incomplete gamma function ratios  \thanks{
This work was supported by  {\emph{Ministerio de Ciencia e Innovaci\'on}}, 
project MTM2009-11686.
 NMT acknowledges financial support from the \textit{Gobierno de Arag\'{o}n}.
}}

\author{Amparo Gil\thanks{
Departamento de Matem\'atica Aplicada y CC. de la Computaci\'on.
ETSI Caminos. Universidad de Cantabria. 39005-Santander, Spain.
({\tt amparo.gil@unican.es}).}
 \and Javier Segura \thanks{Departamento de Matem\'aticas, Estad\'{\i}stica y Computaci\'on.
Facultad de Ciencias. Universidad de Cantabria. 39005-Santander, Spain.
({\tt javier.segura@unican.es}).}
        \and Nico M. Temme\thanks{CWI,
    Science Park 123, 1098 XG Amsterdam, The Netherlands.
({\tt nico.temme@cwi.nl}).} 
}

\begin{document}

\maketitle

\begin{abstract}
Algorithms for the numerical evaluation of
 the incomplete gamma function ratios $P(a,x)=\gamma(a,x)/\Gamma(a)$ and $Q(a,x)=\Gamma(a,x)/\Gamma(a)$

 are described for positive values of $a$ and $x$. Also, inversion methods are 
given for solving the equations $P(a,x)=p$, $Q(a,x)=q$, with $0<p,q<1$. 
Both the direct computation and the inversion of the incomplete gamma function ratios are used
in many problems in statistics and applied probability. The analytical approach from earlier literature is summarized and new initial estimates are derived for starting the inversion algorithms. The performance of the associated software to our algorithms (the Fortran 90 module {\bf IncgamFI}) 
is analyzed and compared with earlier published algorithms. 
\end{abstract}

\begin{keywords} 
Incomplete gamma function ratios, chi-squared distribution function,
inversion of incomplete gamma functions,
numerical evaluation of special functions, 
asymptotic analysis.
\end{keywords}


\pagestyle{myheadings}
\thispagestyle{plain}

\section{Introduction}\label{introd}

As it is well known, the chi-squared distribution (or its equivalent, the
incomplete gamma integral) plays a key role in many applied
probability problems. The incomplete gamma functions are defined by

\begin{equation}\label{eq:int01}
\gamma(a,x)=\int_0^x t^{a-1} e^{-t}\,dt, \quad
\Gamma(a,x)=\int_x^{\infty} t^{a-1} e^{-t}\,dt,
\end{equation}
with ratios
\begin{equation}\label{eq:int02}
P(a,x)=\frac{1}{\Gamma(a)}\gamma(a,x), \quad
Q(a,x)=\frac{1}{\Gamma(a)}\Gamma(a,x),
\end{equation}
where we assume that $a$ and $x$ are positive.

The ratios $P(a,x), Q(a,x)$  are the standard chi-squared probability functions $P(\chi^2| \nu)$ and
$Q(\chi^2 | \nu)$ with parameters $a=\nu/2$ and $x=\chi^2/2$.

Not only the direct computation but also the inversion of cumulative distribution functions is an important
topic in statistics, probability theory, communication theory and econometrics, in
particular for computing percentage points of  the gamma and beta distributions, which have
the chi-square, $F$, and
Student's $t$-distributions as special forms. In the tails of these distributions the
numerical inversion is not very easy, and for these standard distributions
asymptotic formulas are available. 
 
As an example of application in Telecom Engineering, consider a communication system with one
transmit and $r$ receive antennas operating over a flat Rayleigh fading
MIMO channel. For a given communication rate $\Phi$, the outage
probability can be expressed as \cite{Telatar:1999:CMGC}

\begin{equation}\label{P1}
P^{(1\times r)}_{out}=\Frac{1}{\Gamma (r)}\gamma\left(r, \Frac{2^{\Phi}-1}{S/N}\right)
\end{equation}
where $S$ and $N$ are the signal and noise power at the detection moment,
respectively, and $\Gamma$ and $\gamma$ are the gamma  and incomplete gamma functions.

On the other hand, for a communication system with $t$ transmit and one receive antennas
one has that the outage probability is given by \cite{Telatar:1999:CMGC}:

\begin{equation}\label{P2}
P^{(t\times 1)}_{out}=\Frac{1}{\Gamma (t)}\gamma\left(t, \Frac{2^{\Phi}-1}{S/(tN)}\right)\,.
\end{equation}

Then, in order to express the communication rate $\Phi$ in terms of the
desired outage probability $P_{out}$, one has to invert the incomplete gamma function.

In this paper, we present numerical algorithms for the 
computation and inversion of the incomplete gamma function ratios.
A Fortran 90 version of the algorithms is made available at our 
website \footnote{http://personales.unican.es/gila/incgam.zip}.
  
In the numerical algorithms for the
computation of the ratios, both $P(a,x)$ and $Q(a,x)$ are computed. First the primary function (the smaller of the two) is computed, and next the other one by using
\begin{equation}\label{eq:int03}
P(a,x)+Q(a,x)=1.
\end{equation}
In particular for large values of $a,x$ we have a transition at $a\sim x$, with 
\begin{equation}\label{eq:int04}
\begin{array}{ll}
P(a,x)\lesssim\frac12  \quad {\rm when} \quad a\gtrsim x,\\[8pt]
Q(a,x)\lesssim\frac12   \quad {\rm when} \quad a\lesssim x.
\end{array} 
\end{equation}
In the next section we use more refined relations for small values of the parameters.

The algorithms are partly based on \cite{Gautschi:1979:CPI}, where also negative values of $a$ are considered for the  the pair $\{\Gamma(a,z),x^{-a}P(a,x)\}$, the second element being analytic at $x=0$. For applications in probability theory and mathematical statistics we prefer working with the ratios. Also, when $a\sim x$ (both large) the second element becomes approximately $\tfrac12x^{-a}$, and underflow may occur, in which case  $P(a,x)$ cannot be computed. The element $\Gamma(a,x)$ may also become too large to handle. 

For the large parameter case, in particular when $a\sim x$, we do not use Gautschi's approach (continued fractions), but use the method given in \cite{Temme:1987:OTC}; see also \cite[\S8.3]{Gil:2007:NSF} and \cite[\S5.2]{Temme:2007:NAS}. This method is based on uniform asymptotic expansions of the incomplete gamma functions, see \cite{Temme:1979:AEI}. In \cite{Didonato:1986:CIG} these uniform expansions are also used. These authors use several expansions of coefficients for the critical region $a\sim x$, and we will explain that a more efficient expansion can be used.

For the numerical inversion we consider
 the equations
\begin{equation}\label{eq:int05}
P(a,x)=p,\quad Q(a,x)=q, \quad 0<p,\ q<1,
\end{equation}
for a given value of $a$ and we compute $x$. For several cases we give new initial estimates for starting a numerical inversion process. The approach for large $a$ is based on asymptotic methods developed in \cite{Temme:1992:AIG}; see also \cite[\S10.3.1]{Gil:2007:NSF} and \cite[\S6]{Temme:2007:NAS}.

The resulting algorithms and software  significantly improve both the accuracy and ranges of computation of
the algorithm presented in \cite{Didonato:1986:CIG}.  

\section{Methods of computation}\label{methods}
We describe the methods and the domains in the $(x,a)$ quarter plane where they are used. First we define a function
for separating the $(a,x)$ quarter plane in two parts for assigning the primary function, that is, the function $P(a,x)$ or $Q(a,x)$ that has to be computed first.

As a minor modification of the function introduced in \cite{Gautschi:1979:CPI} we use
\begin{equation}\label{eq:moc01}
\alpha(x)= 
\begin{cases} x & \text{if \ $x\ge\tfrac12$,}
\\
\dsp{\frac{\ln\frac12}{\ln (\frac12x)}} &\text{if \ $0<x< \tfrac12$.}
\end{cases}
\end{equation}
Then, the primary function is 
\begin{equation}\label{eq:moc02}
\begin{array}{ll}
P(a,x)   \quad {\rm when} \quad a \ge \alpha(x),\\[8pt]
Q(a,x) \quad {\rm when} \quad a< \alpha(x).
\end{array}
\end{equation}

The function $\alpha(x)$ has the same asymptotic behavior as the one given in \cite{Gautschi:1979:CPI} for small $x$. The
changeover point $x=\tfrac12$ is slightly better.
At $x=a=\tfrac14$ (changeover point in \cite{Gautschi:1979:CPI}) we have $P(a,x)\simeq 0.74$, while at 
the point $x=a=\tfrac12$ we have $P(a,x)\simeq 0.68$.  

In many representations we have to deal with the function
\begin{equation}\label{eq:moc03}
D(a,x)=\frac{x^a e^{-x}}{\Gamma(a+1)},
\end{equation}
see for example representation \eqref{eq:pmoc01}. For large values of $a$ straightforward computation of $\Gamma(a+1)$ will cause overflow. When $x$ is also large, say $x\sim a$, $D(a,x)$ may be computable, as can be seen when using Stirling's formula for the gamma function.

We write $D(a,x)$ in the form
\begin{equation}\label{eq:moc04}
D(a,x)=\frac{e^{-\frac12a\eta^2}}{\sqrt{2\pi a}\, \Gamma^*(a)},
\end{equation}
where
\begin{equation}\label{eq:moc05}
\Gamma^*(a)=\frac{\Gamma(a)}{\sqrt{2\pi/a}\, a^ae^{-a}},\quad a>0,
\end{equation}
and 
\begin{equation}\label{eq:moc06}
\tfrac12\eta^2=\lambda -1-\ln \lambda,\quad \lambda=x/a.
\end{equation}

The quantity $\eta$ arises  in other cases in this paper, and when we take the square root in \eqref{eq:moc06} we assume that for 
$\lambda>0$: $\sign(\eta)=\sign(\lambda-1)$. Then, we have $\eta\sim\lambda-1$ if  $\lambda\sim 1$.

The function $\Gamma^*(a)$ has the asymptotic expansion (Stirling series)
\begin{equation}\label{eq:moc07}
\Gamma^*(a)\sim 1+\tfrac{1}{12}a^{-1}+ \tfrac{1}{288}a^{-2}+\ldots, \quad a \to \infty.
\end{equation}
This function is included in our software package.

For testing the algorithms it is useful to use the recurrence relations
\begin{equation}\label{eq:moc08}
P(a+1,x)=P(a,x)-D(a,x),\quad Q(a+1,x)=Q(a,x)+D(a,x).
\end{equation}
For large values of $a$ and $x$ we can use a scaled version by writing
\begin{equation}\label{eq:moc09}
p(a,x)=P(a,x)/D(a,x),\quad q(a,x)=Q(a,x)/D(a,x), \quad x>0,
\end{equation}
and these functions satisfy the recursion
\begin{equation}\label{eq:moc10}
\frac{x}{a+1}p(a+1,x)=p(a,x)-1,\quad \frac{x}{a+1}q(a+1,x)= q(a,x)+1.
\end{equation}

\subsection{Domains of computation}\label{domains}

\begin{figure}
\begin{center}
\epsfxsize=8cm \epsfbox{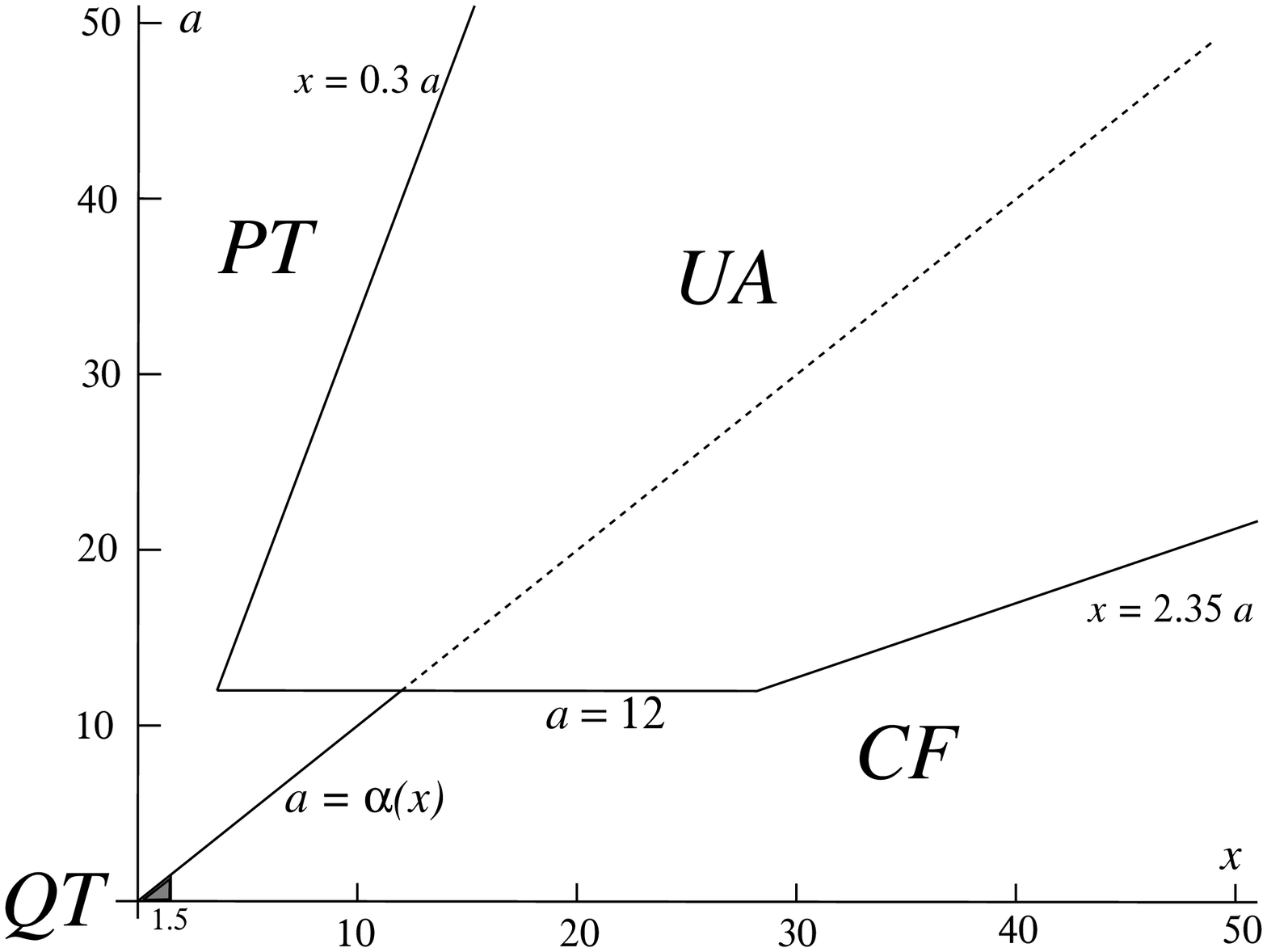}
\caption{The four  domains for computing $P(a,x)$ and $Q(a,x)$.
\label{fig:fig01}}
\end{center}
\end{figure}

In Figure~\ref{fig:fig01} we indicate four domains of computation for the incomplete
gamma function ratios corresponding to different methods: 

\begin{description}
\item[PT:] the Taylor expansion of $P(a,x)$, see \S\ref{ptaylor}.
\item[QT:] the Taylor expansion of $Q(a,x)$, see \S\ref{qtaylor}.
\item[CF:] the continued fraction for $Q(a,x)$, see \S\ref{qconfr}.
\item[UA:] uniform asymptotic methods for  $P(a,x)$ and $Q(a,x)$, see \S\ref{pquni}.
\end{description}

The domains of computation are established following a compromise between efficiency and accuracy: when two methods
provide the same accuracy in a certain parameter region, the selection of one method or another will
depend on the efficiency of each of these methods. The recurrence relations (\ref{eq:moc08}) and (\ref{eq:moc10})
will provide numerical checks for testing the accuracy of the resulting algorithm in all regions of the $(x,a)$-plane.   

\subsection{\protectbold{P(a,x)}: Taylor expansion}\label{ptaylor}

The domain of computation of this section is indicated by {\bf PT} in Figure~\ref{fig:fig01}. 

The expansion is
\begin{equation}\label{eq:pmoc01}
P(a,x)=\frac{x^ae^{-x}}{\Gamma(a+1)}\sum_{n=0}^\infty\frac{x^n}{\left(a+1\right)_n},
\end{equation}
where we use the Pochhammer symbol $\left(a+1\right)_n=\Gamma(a+n+1)/\Gamma(a+1)$.

The series converges for all $a$ and $x$, and the rate of convergence improves as $a/x\to\infty$. The terms of the series are decreasing, because we apply this expansion when $a>x$.

To have an idea about the number of terms needed for obtaining a certain error $\varepsilon$ after truncating the series we write
\begin{equation}\label{eq:pmoc02}
\sum_{n=0}^\infty\frac{x^n}{\left(a+1\right)_n}=S_{n_0}(a,x)+R_{n_0}(a,x),
\end{equation}
where
\begin{equation}\label{eq:pmoc03}
S_{n_0}(a,x)=\sum_{n=0}^{n_0-1}\frac{x^n}{\left(a+1\right)_n},\quad 
R_{n_0}(a,x)=\sum_{n=n_0}^\infty\frac{x^n}{\left(a+1\right)_n},
\end{equation}
and we compute the smallest $n=n_0$ that satisfies 
\begin{equation}\label{eq:pmoc04}
\frac{x^n}{\left(a+1\right)_n}\le\varepsilon.
\end{equation}

In Table~\ref{tab:ig01} we show the smallest number $n$ satisfying \eqref{eq:pmoc04} for $\varepsilon=10^{-15}$ and several values of $a$ and $x$. For large $a$ we see that the number of terms  needed ($n$) becomes constant. This can be understood by observing that the left-hand side of \eqref{eq:pmoc04} becomes roughly $\lambda^n$, $\lambda=x/a$, if $a\gg n$. Because for large values of $a$ we use the method of this section for $P(a,x)$ only if $\lambda\le3/10$,  see Figure~\ref{fig:fig01}, it follows that not more than 30 terms are needed for $a\le10000$.

\begin{table}
\caption{The smallest number $n$ satisfying \eqref{eq:pmoc04}, $\varepsilon=10^{-15}$.
\label{tab:ig01}}
\begin{center}
\begin{tabular}{|| r | c | c | c |c |c |c |c |c ||}
\hline
\quad\quad$x$ & $a/10$ & $  2a/10$ & $ 3a/10$ & $ 4a/10$ & $ 5a/10$& $ 6a/10$& $ 7a/10$& $ 8a/10$\\
 $a\quad$   &  &  & & & & & &\\ \hline
10       & 12  &  15 & 18    & 21 & 23  & 26 & 28 &   31 \\    
100     & 14  &  20 & 26    & 32 & 39  & 47 & 57 &   67 \\    
500     & 14  &  21 & 28    & 36 & 46  & 60 & 79 & 106 \\    
1000   & 14  &  21 & 28    & 36 & 48  & 63 & 86 & 122 \\    
5000   & 14  &  21 & 28    & 37 & 49  & 66 & 94 & 145 \\    
10000 & 14  &  21 & 28    & 37 & 49  & 67 & 95 & 149 \\
\hline
\end{tabular}
\end{center}
\end{table}

The remainder $R_{n_0}$ in \eqref{eq:pmoc03} can be written in the form
\begin{equation}\label{eq:pmoc05} 
R_{n_0}(a,x)=\frac{x^{n_0}}{\left(a+1\right)_{n_0}}\left(a+n_0\right)e^x x^{-a-n_0}\int_0^x
t^{a+n_0-1}e^{-t}\,dt,
\end{equation}
which is again an incomplete gamma function. The  function $u(t)=t-(a+n_0-1)\ln t$ is monotonic in $(0,x]$, and we can integrate with respect to $u$, giving
\begin{equation}\label{eq:pmoc06} 
\int_0^x t^{a+n_0-1}e^{-t}\,dt=\int_{u(x)}^\infty e^{-u} f(u)\,du, 
\end{equation}
where
\begin{equation}\label{eq:pmoc07} 
f(u)=-\frac{dt}{du}=\frac{t}{a+n_0-t-1}\le \frac{x}{a+n_0-x-1}.
\end{equation}
This gives
\begin{equation}\label{eq:pmoc08} 
R_{n_0}(a,x)\le\frac{x^{n_0}}{\left(a+1\right)_{n_0}}\frac{a+n_0}{a+n_0-x-1}.
\end{equation}

\subsection{\protectbold{Q(a,x)}: Taylor expansion}\label{qtaylor}
The expansions of this section will be used for 
\begin{equation}\label{eq:qmoc01}
0\le x\le 1.5, \quad 0\le a\le\alpha(x),
\end{equation}
where $\alpha(x)$ is defined in \eqref{eq:moc01}. This domain is indicated by {\bf QT} in Figure~\ref{fig:fig01}. 
For details and discussion we refer to \cite[\S4.1]{Gautschi:1979:CPI}.

For this case we use the expansion
\begin{equation}\label{eq:qmoc02}
P(a,x)=\frac{x^a}{\Gamma(a)}\sum_{n=0}^\infty \frac{(-1)^n x^n}{(a+n) n!}.
\end{equation}
Straightforward use of the relation $Q(a,x)=1-P(a,x)$ should be avoided when $a$ is small, and we write
\begin{equation}\label{eq:qmoc03}
Q(a,x)= u + v,
\end{equation}
where
\begin{equation}\label{eq:qmoc04}
u=1- \frac{1}{\Gamma(1+a)}+\frac{1-x^a}{\Gamma(1+a)}.
\end{equation}

For the first term we have available an algorithm to compute the function $g(a)$ in the representation
\begin{equation}\label{eq:qmoc05}
1- \frac{1}{\Gamma(1+a)}=a(1-a)g(a).
\end{equation}
The second term can be computed by using an expansion of $1-x^a=1-\exp(a\ln(x))$.

The term $v$ in \eqref{eq:qmoc03} is
\begin{equation}\label{eq:qmoc06}
v=\frac{x^{a}}{\Gamma(1+a)}\left(1-\Gamma(1+a)x^{-a}P(a,x)\right),
\end{equation}
in which we can use the series in \eqref{eq:qmoc02} by skipping the term with $n=0$.

\subsection{\protectbold{Q(a,x)}: continued fraction}\label{qconfr}
The domain of computation of this section is indicated by {\bf CF} in Figure~\ref{fig:fig01}. 

The continued fraction for $Q(a,x)$ is in the form \cite{Gautschi:1979:CPI}
\begin{equation}\label{eq:qmoc08}
Q(a,x)=\frac{x^a e^{-x}}{(x+1-a)\Gamma(a)}\left(
\frac1{1+}\,\frac{a_1}{1+}\,
\frac{a_2}{1+}\,\frac{a_3}{1+}\,
\frac{a_4}{1+}\,\ldots\right),
\end{equation}
where
\begin{equation}\label{eq:qmoc09}
a_k=\frac{k(a-k)}{(x+2k-1-a)(x+2k+1-a)}, \quad k\ge1.
\end{equation}
The front factor term $(x+1-a)$ is not causing problems, because we use the continued  fraction for $x\ge a$. When $a=1,2,3,\ldots$ 
the fraction is terminating. 

Several algorithms are available for the numerical evaluation of continued fractions. See, for example, \cite[\S6.6]{Gil:2007:NSF}. Gautschi used a conversion into an infinite series of the form
\begin{equation}\label{eq:qmoc10}
Q(a,x)=\frac{x^a e^{-x}}{(x+1-a)\Gamma(a)}S(a,x),\quad S(a,x)=1+\sum_{k=1}^\infty t_k,
\end{equation}
where (with $\rho_0=0$)
\begin{equation}\label{eq:qmoc11}
t_k=\rho_1\rho_2\cdots \rho_k,\quad
\rho_k=-\frac{a_k\left(1+\rho_{k-1}\right)}{1+a_k\left(1+\rho_{k-1}\right)},\quad  k\ge1.
\end{equation}

In Table~\ref{tab:ig02} we show the smallest number $n$ satisfying $t_n/S(a,x)<10^{-15}$ for several values of  
$x$ and  $a=x/2.35\rho$; for $\rho$ see the table. 
The $n$-values in this table, as well as those of Table~\ref{tab:ig01}, are not a priori computed but
obtained from numerical computations for the examples in the tables.
When $\rho=1$ and  $a\ge12$, the relation $a=x/{2.35}\rho$ corresponds 
with the border between the domains {\bf UA} and {\bf CF} in Figure~\ref{fig:fig01}. For other values of $\rho$ and 
$x\ge1.5$ the half-lines  $a=x/2.35\rho$ are in the domain  {\bf CF}. Near the diagonal $x=a$ convergence becomes rather slow.

\begin{table}
\caption{The smallest number $n$ satisfying $t_n/S(a,x)<10^{-15}$ for the series in \eqref{eq:qmoc10}, for some $x$ and 
$a={x}/{2.35}\rho$, with several values of $\rho$.
\label{tab:ig02}}
\begin{center}
\begin{tabular}{|| r | r | r | r |r |r |r |r |r ||}
\hline
\quad\quad$\rho$ & $0.1$ & $0.2$ & $  0.4$ & $ 0.6$ & $ 0.8$ & $ 1.0$& $ 2.0$& $ 2.34$\\
 $x\quad$   & &  &  & & & & &\\ \hline
1.5       & 58    & 58  & 57  & 56  &  55 & 54 & 51 &  52  \\    
2          & 45    & 45  & 44  & 43  &  42 & 40 & 40 &  32   \\    
10        & 13    & 12  & 12  & 11  &  11 & 11 & 14 &  14  \\    
100      &  5     & 6    & 7    & 8    &   9  & 11 & 27 &  39  \\    
500      &  4     & 5    & 6    & 6    &   7  &  8  & 25 &  68  \\    
1000    &  4     & 5    & 5    & 6    &   6  &  7  & 19 &  84  \\    
5000    &  4     & 4    & 4    & 5    &   5  &  5  & 10 & 133 \\    
10000  &  3     & 4    & 4    & 4    &   4  &  5  &   9 & 154 \\
\hline
\end{tabular}
\end{center}
\end{table}

The continued fraction is an excellent alternative for the asymptotic expansion
\begin{equation}\label{eq:qmoc12}
Q(a,x)\sim \frac{x^{a-1}e^{-x}}{\Gamma(a)} \sum_{n=0}^{\infty}\frac{(-1)^n(1-a)_n}{x^n},\quad x\to\infty.
\end{equation}

\subsection{\protectbold{P(a,x), Q(a,x)}:  uniform asymptotic expansion}\label{pquni}
The domain of computation of this section is indicated by {\bf UA} in Figure~\ref{fig:fig01}. 
For more details on the used method we refer to \cite{Temme:1987:OTC}; see also \cite[\S8.3]{Gil:2007:NSF} and \cite[\S5.2]{Temme:2007:NAS}. 
We summarize the main steps for constructing an algorithm for this case.

We use the representations \cite{Temme:1979:AEI}
\begin{equation}\label{eq:pqmoc01}
\begin{array}{l}
\dsp{Q(a,x)} = \dsp{\tfrac12\ \erfc(\eta\sqrt{{a/2}}) + R_a(\eta)},\\[8pt]
\dsp{P(a,x)} = \dsp{\tfrac12\ \erfc(-\eta\sqrt{{a/2}}) - R_a(\eta)},
\end{array}
\end{equation}
where 
\begin{equation}\label{eq:pqmoc02}
\erfc\,x=\frac{2}{\sqrt{\pi}}\int_x^\infty e^{-t^2}\,dt,
\end{equation}
the complementary error function. The quantity $\eta$ is defined in \eqref{eq:moc06}, with again $\lambda=x/a$. 

For $R_a(\eta)$ we have  
\begin{equation}\label{eq:pqmoc03}
R_a(\eta) = \frac{e^{-\frac12a\eta^2}}{\sqrt{2\pi a}}S_a(\eta),\quad
S_a(\eta) \sim \sum_{n=0}^\infty \frac{C_n(\eta)}{a^n},
\end{equation}
as $a\to\infty$.
Note that the symmetry relation $P(a,x)+Q(a,x)=1$
is preserved in the  representations in \eqref{eq:pqmoc01} because
$\erfc\,z+\erfc(-z)=2$.

Although analytical expressions for the coefficients $C_n(\eta)$ are available, these representations are difficult to evaluate numerically for small values of $\eta$, that is near the transition $x\sim a$. In \cite{Didonato:1986:CIG} power series expansions of the coefficients $C_n(\eta)$ for $n=0,1,2,\ldots,9$ for small values of $\eta$  are used. 

For the present numerical algorithm we use a different approach.
Instead of expanding each coefficient $C_n(\eta)$ we expand the function
$S_a(\eta)$ of \eqref{eq:pqmoc03} in powers of $\eta$. The
coefficients are functions of $a$, and we write
\begin{equation}\label{eq:pqmoc10}
S_a(\eta)=\sum_{n=0}^\infty \alpha_n \eta^n.
\end{equation}

To compute the coefficients $\alpha_n$, we write 
\begin{equation}\label{eq:pqmoc16}
\alpha_n=\frac{\beta_n}{\Gamma^*(a)},\quad n=0,1,2,\ldots,
\end{equation}
where $\Gamma^*(a)$ is defined in \eqref{eq:moc05}. For $\beta_n$ we have the recursion
\begin{equation}\label{eq:pqmoc17}
\beta_n=\frac1a(n+2)\beta_{n+2}+d_{n+1},\quad n=0,1,2,\ldots.
\end{equation}
where the coefficients $d_n$ are defined as the coefficients in the power series
\begin{equation}\label{eq:pqmoc12}
\frac{\eta}{\lambda-1}=\sum_{n=0}^\infty d_n\eta^n.
\end{equation}
This series  and the one in \eqref{eq:pqmoc10} converge  for  $|\eta|<2\sqrt{\pi}$.

The first values of $d_n$ are
\begin{equation}\label{eq:pqmoc14}
d_0=1,\quad d_1=-\tfrac13,\quad d_2=-\tfrac1{12},\quad
d_3=-\tfrac2{135},\quad
d_4=\tfrac1{864},\quad
d_5=\tfrac1{2835}.
\end{equation}

To describe the algorithm, we choose a positive integer $N$, put $\beta_{N+2}=\beta_{N+1}=0$, and
compute the sequence
\begin{equation}\label{eq:pqmoc18}
\beta_{N},\beta_{N-1},\ldots,\beta_1,\beta_0
\end{equation}
from the recurrence relation
\eqref{eq:pqmoc17}. This recursion is stable in the backward direction.

Because
\begin{equation}\label{eq:pqmoc19}
\alpha_1=\frac{a}{\Gamma^*(a)}\left(\Gamma^*(a)-1\right), \quad
\Gamma^*(a)=1+\frac1a\,\beta_1,
\end{equation}
we have
\begin{equation}\label{eq:pqmoc20}
S_a(\eta) \approx \frac{a}{a+\beta_1}\sum_{n=0}^{N} \beta_n \eta^n
\end{equation}
as an approximation for $S_a(\eta)$.

We use the approximation in \eqref{eq:pqmoc20} for computing the incomplete
gamma functions in IEEE double precision for $a\ge12$ and $|\eta|\le 1$. We need the storage
of 25 coefficients $d_n$, and in  the series in \eqref{eq:pqmoc20}. For $a=12$ we need
25 terms; as $a$ increases the convergence in the algorithm improves and we need a fewer number of terms.

The value $\eta=-1$ corresponds to $\lambda=0.30\ldots,$ and the value $\eta=1$
to $\lambda=2.35\ldots.$ In Figure~\ref{fig:fig01} we show the area indicated by {\bf UA} in the $(x,a)$ quarter-plane
where we can apply the algorithm to obtain IEEE double precision.

\section{Inversion methods}\label{invmethods}

Several approaches are available in the (statistical)
literature for computing the inverse of cumulative distribution
functions, where often a first approximation of $x$ is constructed, based on
asymptotic  estimates, but this first approximation may not be reliable. Higher
approximations
can be obtained by numerical inversion techniques, which require
evaluation of the incomplete gamma functions.
This may be rather time consuming, especially when $a$ is large.

In  \cite{Didonato:1986:CIG} the inversion is considered also; we are using different methods 
based on analytic inversion of power series and asymptotic expansions. In this way it is 
clear how the first steps in the inversion method are taken. We use Newton methods when reliable starting values are available.

\begin{figure}
\begin{center}
\epsfxsize=12cm \epsfbox{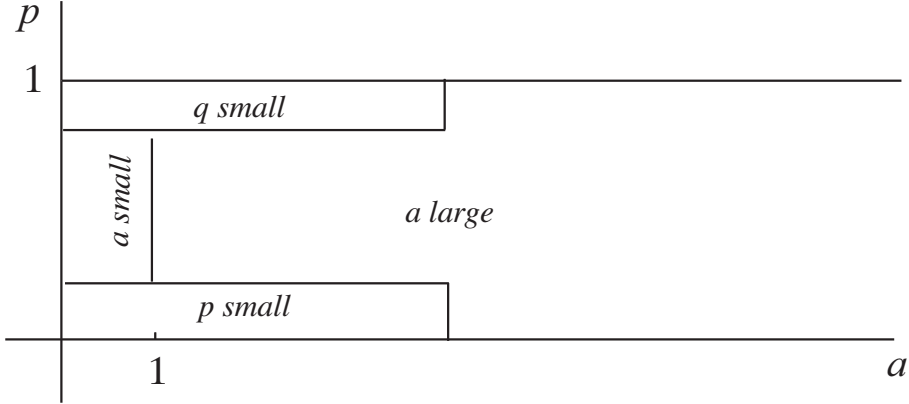}
\caption{The four  domains for inverting $P(a,x)=p$ and $Q(a,x)=q$ with $p+q=1$.
\label{fig:fig02}}
\end{center}
\end{figure}

We solve the equations
\begin{equation}\label{eq:pqinv01}
P(a,x)=p,\quad Q(a,x)=q,\quad 0\le p\le 1, \quad 0\le q\le1,
\end{equation}
for $x$, with $a$ as a given positive
parameter. We consider several cases, which are schematically indicated in 
Figure~\ref{fig:fig02}. In most cases we invert the equation with $\min(p,q)$.
Of course, if $x(p,a)$ denotes the solution of the first equation then the solution of the second equation satisfies
$x(q,a)=x(1-p,a)$. We assume that the user provides both $p$ and $q$, which is important when $\min(p,q)$ is small.

\subsection{Small values of \protectbold{p}}
\label{invpsmall}
When $p$ is small we use the series in \eqref{eq:qmoc02}, and write the inversion problem as
\begin{equation}\label{eq:pqinv02}
x=r\left(1+\sum_{n=1}^\infty\frac{a(-1)^n x^n}{(a+n)n!}\right)^{-1/a},\quad r=\left(p\Gamma(1+a)\right)^{1/a},
\end{equation}
where we assume that $r$ is small. Inverting this relation, we obtain $\dsp{x=r+\sum_{n=2}^\infty c_kr^k}$, and the first few coefficients are
\begin{equation}\label{eq:pqinv03}
    \begin{array}{l}
c_2= \dsp{\frac{1}{a+1},}\\[8pt]
c_3= \dsp{\frac{3a+5}{2(a+1)^2(a+2)},}\\[8pt]
c_4= \dsp{\frac{8a^2+33a+31}{3(a+1)^3(a+2)(a+3)},}\\[8pt]
c_5= \dsp{\frac{125a^4+1179a^3+3971a^2+5661a+2888}{24(1+a)^4(a+2)^2(a+3)(a+4)}.}
\end{array}
\end{equation}

It appears that  $c_k=\bigO((a+1)^{-k+1}$ for large values of $a$, and from numerical experiments 
by checking several values of $p$, $q$, $a$ which were relevant for these cases,
we conclude that if 
$r<0.2(1+a)$, that is, $p<(0.2(1+a))^{a}/\Gamma(1+a)$ we can obtain 4 digits accuracy in $x$ 
with the coefficients shown in \eqref{eq:pqinv03}. This is enough for starting a Newton method for obtaining higher accuracy.

This method also works when $a$ is small, because in that case $r\sim\exp(\frac1a\ln p)$ becomes small for all fixed $p\in(0,1)$.

\subsection{Small values of \protectbold{q}}\label{invqsmall}
When $q$ is small we use the asymptotic expansion in \eqref{eq:qmoc12}. A first approximation $x_0$  of $x$ is obtained from the equation
\begin{equation}\label{eq:pqinv04}
e^{-x_0}x_0^{a}=q\Gamma(a)
\end{equation}
Higher approximations of $x$ are obtained in the form $\dsp{x\sim x_0-L+b\sum_{k=1}^\infty d_k/x_0^k}$, where  $b=1-a$, $L=\ln(x_0)$, with first coefficients
\begin{equation}\label{eq:pqinv05}
    \begin{array}{l}
d_1= L-1,\\[8pt]
d_2= \tfrac12 (3b-2bL+L^2-2L+2),\\[8pt]
d_3= \tfrac1{6}(24bL-11b^2-24b-6L^2+ \\[8pt]
\quad\quad\quad12L-12-9bL^2+6b^2L+2L^3),\\[8pt]
d_4= \tfrac1{12}(72+36L^2+3L^4-72L+162b-168bL-12L^3+\\[8pt]
\quad\quad\quad
25b^3-22bL^3+36b^2L^2-12b^3L+84bL^2+120b^2-114b^2L).
\end{array}
\end{equation}

These coefficients, as well as those given in (\ref{eq:pqinv03}), are obtained by symbolic computation.
Maple codes for computing these or more coefficients can be obtained
at our website \footnote{http://personales.unican.es/gila/coefMaple.zip}.

This method works with rather small values of $q$ (large values of $x_0$). When we assume that $0<a<10$ and  $q<e^{-\frac12a}/\Gamma(a+1)$
 we obtain about 4 digits accuracy in $x$, which  is enough for starting a Newton method for obtaining higher accuracy. For larger 
values of $a$ and $q$, the method of \S\ref{invalarge} can be used.
\subsection{Small values of \protectbold{a}}\label{invasmall}
We consider the inversion of $P(a,x)=p$ for $a\in(0,1)$. For these values of $a$ we cannot derive expansions for obtaining a reliable starting value, unless $p$ is small in which case we can use the result of \S\ref{invpsmall}. Also, the method of the next section cannot be used in this case.

We observe that
\begin{equation}\label{eq:pqinv06}
P(a,x)=\frac1{\Gamma(a)}\int_0^x t^{a-1}e^{-t}\,dt < \frac1{\Gamma(a)}\int_0^x t^{a-1}\,dt =
\frac{x^a}{\Gamma(a+1)}
\end{equation}
and 
\begin{equation}\label{eq:pqinv07}
P(a,x)>\frac1{\Gamma(a)}\int_0^x e^{-t}\,dt = \frac1{\Gamma(a)}\left(1-e^{-x}\right).
\end{equation}

Let $x_l, x_u$ be defined by 
\begin{equation}\label{eq:pqinv08}
x_l=\left(p\Gamma(a+1)\right)^{1/a},\quad x_u=-\ln\left(1-p\Gamma(a+1)\right).
\end{equation}
Then the solution $x$ of $P(a,x)=p$ with $a\in(0,1)$ satisfies $x_l<x<x_u$. The same results hold for the inversion of $Q(a,x)=q$ for $a\in(0,1)$ when $p$ is replaced with $1-q$.   

These bounds of $x$ can be used for starting values for the Newton method. Of the two possibilities, the best
option is $x_l$ because the Newton method necessarily converges from this starting value. The reason is that $P(a,x)$ is an increasing function
with negative second derivative; elementary graphical arguments show that if an starting value $x_0<x$ is chosen, then the Newton iteration
$x_{n+1}=x_n-P(a,x_n)/P'(a,x_n )$ produces a monotonically increasing sequence which is bounded by $x$, and therefore a converging sequence. 
The same is true for $Q(a,x)$ because it is decreasing and with positive second derivative.

For this case, the initial approximation may be inaccurate, however convergence is certain and not very expensive, as numerical experiments show.

\subsection{Large  values of \protectbold{a}}\label{invalarge}

We perform the inversion of \eqref{eq:pqinv01} with respect to the parameter
$\eta$ by using the representations \eqref{eq:pqmoc01}. Afterwards we
have to compute $\lambda$ and $x$ from the relation for $\eta$ in \eqref{eq:moc06}
and $\lambda=x/a$. We concentrate on the second equation
in \eqref{eq:pqinv01}. For details  we refer to \cite{Temme:1992:AIG}; see also \cite[\S10.3.1]{Gil:2007:NSF} and \cite[\S6]{Temme:2007:NAS}.

We rewrite the inversion problem in the form
\begin{equation}\label{eq:pqinv09}
\tfrac12\erfc\left(\eta\sqrt{{a/2}}\right)+R_a(\eta)=q, \quad q\in[0,1],
\end{equation}
which is equivalent to the second equation in \eqref{eq:pqinv01}, and we denote the
solution of the above equation by $\eta(q,a)$.

To start the procedure we consider $R_a(\eta)$ in  \eqref{eq:pqinv09} as a
perturbation, and we define the number $\eta_0=\eta_0(q,a)$ as the real number
that satisfies the equation
\begin{equation}\label{eq:pqinv10}
\tfrac12\erfc\left(\eta_0\sqrt{{a/2}}\right)=q.
\end{equation}
Computation of $\eta_0$  requires an inversion of the complementary error
function, which is discussed in \S\ref{inverfc}.

For large values of $a$ the value $\eta$ defined by (\ref{eq:pqinv09}) can be approximated by
$\eta_0$:
\begin{equation}\label{eq:pqinv12}
\eta(q,a)=\eta_0(q,a)+\varepsilon(\eta_0,a),
\end{equation}
and  it is possible to expand 
\begin{equation}\label{eq:pqinv13}
\varepsilon(\eta_0,a)\sim\frac{\varepsilon_1(\eta_0,a)}a+\frac{\varepsilon_2(\eta_0,a)}{a^2}+
\frac{\varepsilon_3(\eta_0,a)}{a^3}+\cdots,
\end{equation}
as $a\to\infty$. The coefficients $\varepsilon_j(\eta_0,a)$ can be written explicitly as functions of
$\eta_0$, the first coefficient being
\begin{equation}\label{eq:pqinv14}
\varepsilon_1(\eta_0,a)=\frac1\eta_0\ln \frac{\eta_0}{\lambda_0-1},
\end{equation}
where $\lambda_0$ follows from inverting \eqref{eq:moc06} with $\eta$ replaced by $\eta_0$.

From numerical tests  it follows that we can obtain 3 or 4 significant digits when using \eqref{eq:pqinv13} with 4 terms  for $a\ge 1$ and $q\in(0,1)$. This is enough to start a stable Newton method.

\begin{remark}\label{rem01}
{\rm
We start the inversion of $P(a,x)=p$ with the equation
\begin{equation}\label{eq:pqinv15}
\tfrac12\erfc\left(-\eta\sqrt{{a/2}}\right)-R_a(\eta)=p, \quad p\in[0,1],
\end{equation}
and with the first approximation  $\eta_0(p,a)=-\eta_0(q,a)$ (cf.~\eqref{eq:pqinv10}),
with $p+q=1$. The results for $P(a,x)=p$ then follow from the  results of this section with 
$\eta_0(q,a)$ replaced by $-\eta_0(p,a)$, throughout.
}
\end{remark}

\subsubsection{Complementary error function}
\label{inverfc}
The inversion of this function is the first step 
in the large $a$ asymptotic inversion method for the incomplete gamma function ratios.
We summarize results from \cite[\S7.17]{Temme:2010:ERF}; for more details, see \cite[\S10.2]{Gil:2007:NSF}.

We denote the inverse of the function $x=\erfc\,y$ by $y=\inverfc\,x$. Then, with $t=\frac12\sqrt\pi\,(1-x)$, we have
\begin{equation}\label{eq:inverfc1}
\inverfc\,x=t+\tfrac13t^3+\tfrac{7}{30}t^5+\tfrac{127}{630}t^7+\ldots, \quad 0<x<2.
\end{equation}
This expansion is actually for the inversion of $\erf\,y=x$, useful for small values of $t$. More coefficients in the expansion can be found in \cite{Strecok:1968:OCI}, but they can easily be obtained by using computer algebra and formal manipulation of power series.

For small values of $x$ we have an asymptotic expansion. 
Let $t$, $\alpha$, and $\beta$ be defined by
\begin{equation}\label{eq:inverfc2}
t=\frac{2}{\pi x^{2}},\quad \alpha=\frac1{\ln t}, \quad \beta =
\ln(\ln t).
\end{equation}
Then we have the expansion
\begin{equation}\label{eq:inverfc3}
\inverfc\,x\sim
\frac1{\sqrt{2\alpha}}\left(1+x_{1}\alpha+x_2\alpha^2+x_3\alpha^3+x_4\alpha^4+\cdots\right).
\end{equation}
The first coefficients $x_k$ are given by
\begin{equation}\label{eq:inverfc4}
    \begin{array}{l}
x_1= \dsp{-\tfrac12\beta,}\\[8pt]
x_2= \dsp{-\tfrac18\left(\beta^2-4\beta+8\right),}\\[8pt]
x_3= \dsp{ -\tfrac1{16} \left(\beta^{3}-8\beta^2+32 \beta-56\right),}\\[8pt]
x_4=  \dsp{-\tfrac1{384} \left(15\beta^{4} -184 \beta^3+1152
\beta^2-4128 \beta+7040\right).}
\end{array}
\end{equation}

In \S\ref{sec:cerfh} we use a higher-order Newton process to obtain better approximations.

For other methods of the inversion of the error functions,
we refer to \cite{Strecok:1968:OCI}, where coefficients of the
Maclaurin expansion of $y=\inverf\, x$, the
inverse of $x=\erf\,y$, are given, with Chebyshev coefficients for an
expansion on the $y$-interval $[-0.8,0.8]$.
For small values of $x$
(not smaller than $10^{-300}$) high-precision coefficients of Chebyshev
expansions
are given for the numerical evaluation of  $y=\inverfc\,x$.
For rational Chebyshev (near-minimax) approximations for  $y=\inverfc\,x$,
we refer to \cite{Blair:1976:RCA}, where $x$-values are considered
in the $x$-interval $[10^{-10000},1]$, with relative errors ranging
down to $10^{-23}$. An asymptotic formula for the region
$x\rightarrow 0$ is also given.

\section{High order Newton-like methods}\label{newtonhigh}
\index{Newton-Raphson method!high-order inversion}
We describe the known method (see \cite[\S10.6]{Gil:2007:NSF} and \cite[\S6.8]{Temme:2007:NAS}) for constructing Newton-like methods of high order, with details for the inversion of the complementary error function and the incomplete gamma functions. For certain values of the parameters we will use these methods in the inversion algorithms. 

The method is based on inverting the power series in $h$:\begin{equation}\label{eq:newhigh1}
f(\zeta)=f(\zeta_0+h)=f(\zeta_0)+h f_1+\frac1{2!}h^2 f_2+\frac1{3!} f_3+\cdots\,,
\end{equation}
where $\zeta$ is a zero of $f$ and  $\zeta_0$  is an approximation. Also,  $f_k=f^{(k)}(\zeta_0)$.
We expand, assuming $f(\zeta_0)$ is small, 
\begin{equation}\label{eq:newhigh2}
h=c_1f(\zeta_0)+c_2 f^{2}(\zeta_0)+c_3f^{3}(\zeta_0)+\cdots.
\end{equation}
and find, when $f_1\ne0$,
\begin{equation}\label{eq:newhigh3}
\begin{array}{ll}
\dsp{c_1=-\frac{1}{f_1}},
&
\dsp{c_2=-\frac{f_2}{2f_1^3}},\\[10pt]
\dsp{c_3=\frac{-3f_2^2+f_3^3f_1}{6f_1^5},}
&
\dsp{c_4=-\frac{f_4f_1^2+15 f_2^3-10 f_2 f_3 f_1}{24f_1^7}}.
\end{array}
\end{equation}
When we neglect in \eqref{eq:newhigh2} the coefficients $c_k$ with $k\ge2$, we obtain Newton's Rule,
with $\zeta\dot=\zeta_0-f(\zeta_0)/f^{\prime}(\zeta_0)$.

When $f(z)$ satisfies a simple  ordinary differential equation, the higher derivatives can be replaced
 by combinations of lower derivatives.

\subsection{Complementary error function}\label{sec:cerfh}
We compute a zero of the function $f(y)=\erfc\,y-x$ with $0<x<2$. Assume we have a first approximation $y_0$ then the zero $y$ is written as $y=y_0+h$, where $h$ is as in \eqref{eq:newhigh2} with $f(\zeta_0)=f(y_0)$. The derivatives of $f(y)$ are in terms of Hermite polynomials\footnote{http://dlmf.nist.gov/18.5.E5}:
\begin{equation}\label{eq:newhigh4}
f_k=c \,e^{-y^2}(-1)^{k-1}H_{k-1}(y), \quad k=1,2,3,\ldots,\quad c=-\frac{2}{\sqrt{\pi}}.
\end{equation}
This gives 
\begin{equation}\label{eq:newhigh5}
\begin{array}{l}
c_1=-1/f_1={-e^{y_0^2}/c},\quad
c_2={y_0c_1^2}, \quad 
c_3={\frac13(4y_0^2+1)c_1^3},\\[8pt]
c_4={\frac16y_0(7+12y_0^2)c_1^4}, \quad
c_5=\frac1{30}(8y_0^2+7)(12y_0^2+1)c_1^5.
\end{array}
\end{equation}

\subsection{Incomplete gamma function}\label{sec:inchigh}

When starting with an initial value $x_0>0$ for the inversion of $f(x)=P(a,x)-p$, the coefficients $c_k$ can be derived from the integral representation of $P(a,x)$. This gives
\begin{equation}\label{eq:newhigh6}
\begin{array}{l}
c_1=\dsp{-x_0^{1-a}e^{x_0}\Gamma(a)},\quad c_2=\dsp{\frac{x_0+1-a}{2x_0}c_1^2},\\[6pt]
c_3=\dsp{\frac{2x_0^2+(a-1)(2a-1-4x_0)}{6x_0^2}c_1^3},\\[6pt]
c_4=\dsp{\frac{6x_0^3-(a-1)(18x_0^2-x_0(18a-11)+(2a-1)(3a-1))
}{24x_0^3}c_1^4}.
\end{array}
\end{equation}

For  $a=\frac12$, in which case $P(\frac12,x)=\erf\,\sqrt{x}$, and $p=\frac12$, we take 
$x_0$ of \eqref{eq:pqinv08}, giving $x_0=\pi/16=0.19634954$ and
$f(x_0)=\erf\,\sqrt{x_0}-\frac12=-0.03088405$. Using 4  terms $c_k$
 from \eqref{eq:newhigh6},
it follows that $h= 0.03111855$, giving $x\,\dot=\,x_0+h=0.227468092$,
and with this value we have $f(x)\dot= -1.12\,10^{-7}$.  

For much smaller values of $a$, this fourth order method has some convergence problems. However, as commented before, the plain Newton
method converges with certainty. For values of $a$ smaller than 0.05 the plain Newton is used in the algorithm. 

\section{Associated Software and Testing}\label{numtest}

A Fortran 90 module ({\bf IncgamFI}) implementing the algorithms can be obtained at our 
website \footnote{http://personales.unican.es/gila/incgam.zip}.
The module includes the public routines {\bf incgam}, for the computation of $P(a,x)$ and $Q(a,x)$, and {\bf invincgam},
for the computation of $x$ in the equations $P(a,x)=p$ and $Q(a,x)=q$ with $a$ as a given positive parameter;
$p$ and $q$ are also inputs of this routine.
As a test of the routine {\bf incgam}, we compute the maximum relative errors 
for the relations given in equation  (\ref{eq:moc08}) using $10^6$ and $10^7$ random points
in the following two regions of the $(x,a)$-plane, respectively. We obtain:

\begin{enumerate}
\item{$(0,1]\times (0,1]$:} $ 1.7\, 10^{-15}$,
 \item{$(0,500]\times (0,500]$:} $ 7.9\, 10^{-13}$.
\end{enumerate}

These errors constitute the accuracy claim of our algorithm {\bf incgam}. 

Complementarily, we have performed tests in larger ranges by using the scaled expressions given in (\ref{eq:moc09}) 
and considering the relation (\ref{eq:moc10})
in the regions where the function $D(a,x)$ (\ref{eq:moc03}) can be explicitly factored out. 
The use of scaled functions allows to test these methods for large values of $x$ and $a$. 
We consider first the cases of the Taylor expansion 
 and the continued fraction method.
Using  $10^7$ and $10^8$ random points
in the following two regions of the $(x,a)$-plane (excluding the points where asymptotic expansions are used, see Figure~\ref{fig:fig01}), we obtain as maximum relative errors:

\begin{enumerate}
\item{$(0,10^4]\times (0,10^4]$:} $ 8.3\, 10^{-15}$,
 \item{$(0,10^5]\times (0,10^5]$:} $ 9.1\, 10^{-15}$.
\end{enumerate}

For testing the uniform asymptotic expansions for large values of the parameters,
in the numerical algorithm we use the representations
\begin{equation}\label{eq:remc01}
\begin{array}{l}
\dsp{Q(a,x)} = \dsp{\Gamma^*(a)D(a,x)\left(\tfrac12\sqrt{2\pi a}\ \Erfc(\eta\sqrt{{a/2}}) + S_a(\eta)\right)},\\[8pt]
\dsp{P(a,x)} =  \dsp{\Gamma^*(a)D(a,x)\left(\tfrac12\sqrt{2\pi a}\ \Erfc(-\eta\sqrt{{a/2}}) - S_a(\eta)\right)},
\end{array}
\end{equation}
where $\Erfc\,x=e^{x^2}\erfc\,x$, $D(a,x)$ is defined in \eqref{eq:moc03}, and $S_a(\eta)$ denotes the series as in \eqref{eq:pqmoc03}.
In this way we can test the algorithm for the scaled function $q(a,x)=Q(a,x)/D(a,x)$ by using the relation in \eqref{eq:moc10}. Similarly
for $p(a,x)=Q(a,x)/D(a,x)$.  Using  $10^7$ random points in the region $(0,10^4]\times (0,10^4]$ of the $(x,a)$-plane, the maximum relative
error obtained when computing the scaled functions using uniform asymptotic expansions is $4\, 10^{-14}$.

The results obtained with the scaled expressions confirm the stability and accuracy of
the methods used for computing $P(a,x)$ and $Q(a,x)$.

\begin{figure}
\begin{center}
\epsfxsize=12cm \epsfbox{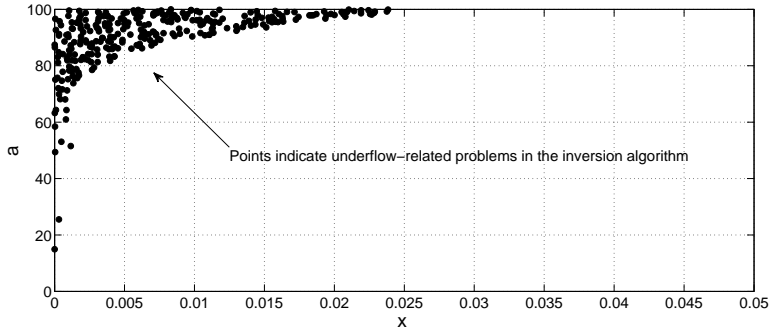}
\caption{Test of the inversion algorithm: underflow-related problems in double precision
arithmetic in the region $(x,a) \in (0,100) \times (0,100)$. Notice that the problems
are concentrated in the interval $0<x<0.025$. 
\label{fig:fig04}}
\end{center}
\end{figure}

For the inversion algorithm, testing is made by checking that the composition of the functions with their inverse is the identity:
from the values of $a$ and $x=x_{in}$ we compute $p=P(a,x_{in})$ and $q=Q(a,x_{in})$ and from the values of $p$, $q$ and $a$ we obtain
$x_{out}$ from the inversion algorithm; we compare the values $x_{in}$ and $x_{out}$ and compute the relative errors.
As a first test, we generate $10^7$ random points in the region $(x,a) \in (0,100]\times (0,100]$. In the
test we exclude the points
where problems related to the underflow limit in double precision arithmetic when computing
the incomplete gamma function ratios appear (see Figure \ref{fig:fig04}). 
The maximum relative error is $1.42\, 10^{-11}$. It is important to note that this value is obtained for a point $(x,a)$ 
near the region of underflow
problems: the value of the function $P(a,x)$ at that point was $\sim 10^{-297}$. Some loss of accuracy is expected
in these cases.

\begin{figure}
\begin{center}
\epsfxsize=12cm \epsfbox{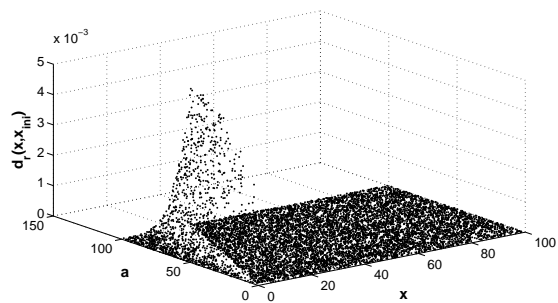}
\caption{Accuracy of the initial estimates for the inversion problem: points correspond to relative distances 
$d_r(x,x_{ini})=\left\|1-x_{ini}/x\right\|$
between the initial estimates ($x_{ini}$) and the true values $x$.
\label{fig:fig05}}
\end{center}
\end{figure}

The accuracy of the initial estimates discussed in sections (\ref{invpsmall}), (\ref{invqsmall}) and
(\ref{invalarge}) is illustrated in Figure \ref{fig:fig05}, where relative distances between the initial
estimate ($x_{ini}$) and the true value ($x$) are plotted. $10^4$ random points have been considered in the plane
  $(x,a) \in (0,100) \times (1,100)$. As can be seen, the poorest estimate in the test is located at a relative distance less than $5\, 10^{-3}$
to the real value. 

\begin{figure}
\begin{center}
\epsfxsize=12cm \epsfbox{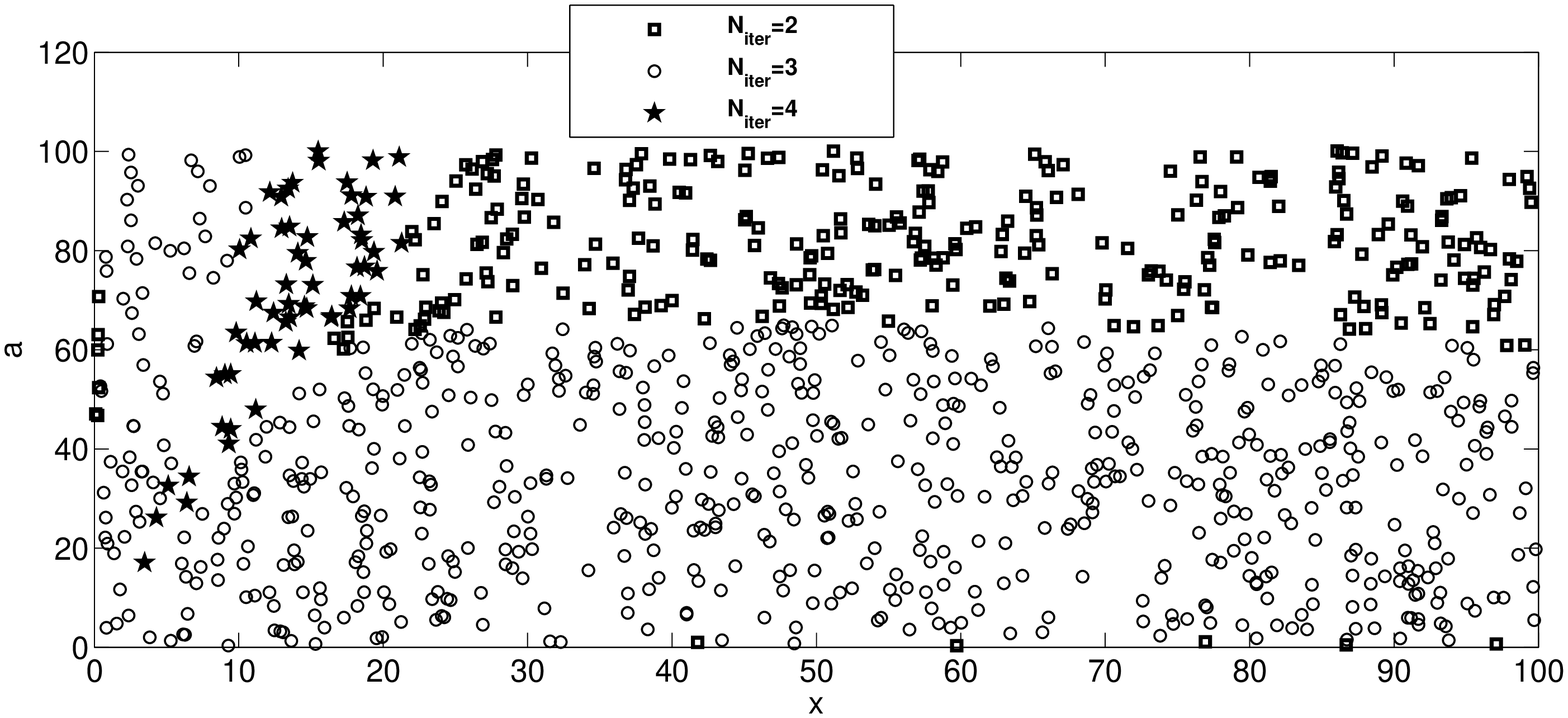}
\caption{Number of iterations used in the inversion algorithm in the region $(x,a) \in (0,100]\times (0,100]$. 
\label{fig:fig06}}
\end{center}
\end{figure}

The number of iterations used in the inversion algorithm is also tested. Figure \ref{fig:fig06}
shows the number of iterations used in the region $(x,a) \in (0,100]\times (0,100]$ for computing 
the $x$ values in the equations $p=P(a,x)$ and $q=Q(a,x)$ within an accuracy of $10^{-12}-10^{-14}$. As can be seen, 2 or 3 iterations are enough
 in most of the points of $(x,a)$-plane.

\begin{figure}
\begin{center}
\epsfxsize=12cm \epsfbox{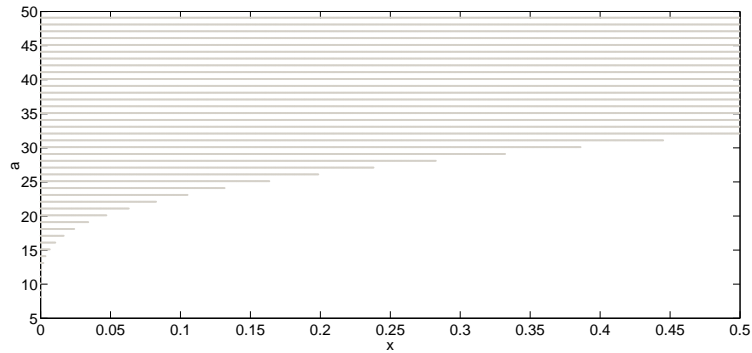}
\caption{Inversion routine {\bf GAMINV} of reference \cite{Didonato:1986:CIG}: points are plotted when the relative accuracy 
in the inversion test is greater than 0.1. The routine
{\bf GRATIO} was used for the direct computation of the functions.
\label{fig:fig07}}
\end{center}
\end{figure}

As a final comment, our module {\bf IncgamFI} clearly improves the algorithm provided in \cite{Didonato:1986:CIG}, having
its associated Fortran 77 routines {\bf GRATIO} and {\bf GAMINV} (single precision routines) both a more
limited range of validity and accuracy than our algorithms.
 As an illustration, Figure \ref{fig:fig07}
shows the points where the inversion routine {\bf GAMINV} fails when performing the same test as used
for our inversion algorithm. A point is plotted when the relative accuracy was greater than 0.1. The routine
{\bf GRATIO} was used for the direct computation of the functions.

\section*{Acknowledgments}
This work was supported by  {\emph{Ministerio de Ciencia e Innovaci\'on}}, 
project MTM2009-11686.

\bibliographystyle{plain}

\end{document}